# The formal roots of Platonism

Bhupinder Singh Anand

We present some arguments for the thesis that a set-theoretic inspired faith, in the ability of intuitive truth to faithfully reflect relationships between elements of a Platonic universe, may be as misplaced as an assumption that such truth cannot be expressed in a constructive, and effectively verifiable, manner.

## 1. An implicit thesis, and an explicit belief

In a 1991 lecture on "The Future of Set Theory", Saharon Shelah presents an overview of classical Set Theory that is based on an implicit thesis - that intuitive truth is essentially non-verifiable - and on the explicit belief that:

> ... ZFC exhausts our intuition except for things like consistency statements, so a proof means a proof in ZFC ... all of us are actually proving theorems in ZFC.[1]

It is not obvious whether Shelah makes a precise distinction between such, vision-based, intuitive truth, and Tarski's formal, but essentially non-constructive, classical definitions of the truth of formulas of a formal system under a given interpretation[2], since he remarks that:

> I am in my heart a card-carrying Platonist seeing before my eyes the universe of sets ... (*regarding*) the role of foundations, and philosophy ... I do not have any objection to those issues per se, but I am suspicious ... My feeling, in an overstated form, is that beauty is for eternity, while philosophical value follows fashion.

---

1 All quotations attributed to Shelah are from [Sh91].

2 Cf. Mendelson ([Me64], p49-53). We take Mendelson [Me64] as a representative exposition of the foundational concepts of classical mathematics considered in this paper.

## 1.1 Shelah's faith may be misplaced

However, Shelah's faith, in the ability of intuitive truth to faithfully reflect relationships between elements of a Platonic universe, may be as misplaced as his assumption that such truth cannot be expressed in a constructive, and effectively verifiable, manner.

As we show in Anand [An02a], there is a constructive, and intuitionistically unobjectionable, proof that the replacement axiom[3] of ZFC is inconsistent with an interpretation of standard PA in ZFC. Thus the question of intuitive truth may be linked to that of the consistent introduction of mathematical concepts into ZFC, through axiomatic postulation, in ways that may not be immediately obvious to a self-confessed Platonist, such as Shelah, even if we grant him the vision that is implicit in his following remarks:

> From the large cardinal point of view: the statements of their existence are semi-axioms, (for extremists - axioms). Adherents will probably say: looking at how the cumulative hierarchy is formed it is silly to stop at stage *omega* after having all the hereditarily finite sets, nor have we stopped with Zermelo set theory, having all ordinals up to *aleph_omega*, so why should we stop at the first inaccessible, the first Mahlo, the first weakly compact, or the first of many measurables? We are continuing the search for the true axioms, which have a strong influence on sets below (even on reals) and they are plausible, semi-axioms at least.

---

[3] We take the replacement (comprehension/separation) axiom (axiom schema) as stating, essentially, that the range of every function, whose domain is a well-defined set in an axiomatic set theory Z, and whose values are always elements of a well-defined set in Z, is a well-defined set in Z.

A very interesting phenomenon, attesting to the naturality of these axioms, is their being linearly ordered (i.e., those which arise naturally), though we get them from various combinatorial principles many of which imitate *aleph*$_0$, and from consistency of various "small" statements. It seems that all "natural" statements are equiconsistent with some large cardinal in this scale; all of this prove their naturality.

This raises the question:

ISSUE: Is there some theorem explaining this, or is our vision just more uniform than we realize?

Intuition tells me that the power set and replacement axioms hold, as well as choice (except in artificial universes), whereas it does not tell me much on the existence of inaccessibles. According to my experience, people sophisticated about mathematics with no knowledge of set theory will accept ZFC when it is presented informally (and well), including choice but not large cardinals. You can use collections of families of sets of functions from the complex field to itself, taking non-emptiness of cartesian products for granted and nobody will notice, nor would an *omega*-fold iteration of the operation of forming the power set disturb anybody. So the existence of a large cardinal is a very natural statement (and an interesting one) and theorems on large cardinals are very interesting as implications, not as theorems (whereas proving you can use less than ZFC does not seem to me very interesting).

Prima facie, Shelah seems to implicitly imply that there is a perceptible difference between treating an assertion as an implication (presumably from a non-intuitive set of axioms and semi-axioms), and treating it as a theorem (presumably from an intuitive set of axioms); however, he also appears to comfortably accept that the dividing line between what is non-intuitive, and what is intuitive, may be essentially subjective.

Now, from within the framework of a Platonist philosophy, differences in individual perceptions of such a dividing line are, perhaps, acceptable as reflecting differences in the levels of awareness, between individual perceivers, of the Platonic existence of mathematical concepts. However, apart from the ethical questions implicit in the acceptance of, firstly, an objective, absolute, truth and, secondly, subjective, gradable, levels of awareness[4], such differences would be a matter of serious concern for disciplines for whom mathematics is, essentially, a language of reliable, and verifiable, external expression and communication. Such a language should, clearly, be based on notions of formal truth that offer greater precision in, and verifiability of, its assertions than that suggested by Shelah's notion of intuitive truth in his above comments.

Perhaps the point is better illustrated by observing how practicing scientists actually view mathematics. For instance, in the pre-print of his forthcoming book, "What is and what will be - Integrating spirituality and science", Budnik [Bu01] writes:

(*a*) If infinity is a potential and never a completed reality then infinite sets do not exist.

(*b*) The objects definable within a formal mathematical system no matter what axioms of infinity it includes are countable (they can be mapped onto the integers). This result is called the Lowenheim-Skolem Theorem. The idea of the proof is that a formal system can be interpreted as a computer program for generating theorems. Such a program can output all of the names of the objects or sets definable with the system. These names and thus the collection of all objects they refer to are countable. They can be mapped onto the integers.

---

[4] As are implicitly implied by subjective expressions such as "people sophisticated about mathematics".

(*c*) All real numbers and for that matter larger cardinals that can ever be defined in any mathematical system that finite creatures create will be countable. They will not be countable from within the system. Cantor`s proof is correct as a proof about formal systems. If real numbers do not exist Cantor's proof is about the structure of formal systems and not some greater metaphysical reality.[5]

(*d*) This suggests that the theory of cardinals is an illusion. It is talking indirectly about ways of extending mathematics that are countable and reducible to properties of computer programs. The set of reals definable within a formal system is a countable set in a more powerful formal system. In the more powerful system there is a countable ordinal that characterizes this set. Because we are talking about names generated by a computer program there is a connection between this countable ordinal and properties of computer programs.

(*e*) It is easy to talk about a formal system plus an uncountable number of axioms that state the existence of all reals. Each real number can be defined as an infinite sequence of digits. We cannot write the entire sequence but we have a sense of what an arbitrary real is. In this way mathematicians talk about the true set theory that includes all reals.

(*f*) The human mathematical mind is the product of biological evolution. There is no evidence of a special facility that transcends the finite. On the contrary all the evidence suggests the opposite. The current 'theological' approach to mathematical truth flies in the face of the evidence. We believe the cardinality of the real numbers is *C*. We do not think any real number as a completed infinite sequence exists. Of course the integers and rational numbers are also considered reals so one could argue that the cardinality of the reals is *aleph*$_0$ or the same as *omega*.

---

[5] Given the present context, it is interesting to assess the meaning of the last two remarks - and the one following - in view of the arguments in Anand [An03b].

(*g*) Mathematics involves the creation of truth which has an objective meaning. This is truth about what a computer does if it is allowed to run forever perhaps following an ever expanding number of paths.

(*h*) Mathematical statements that cannot be interpreted as questions about events all of which will occur in a potentially infinite deterministic universe are neither true nor false in any absolute sense. Of course, they may be useful properties that are either true or false relative to a particular formal system.

(*i*) Mathematics gives us some sense of what is possible and establishes some of the conditions necessary to realize those possibilities.

What is interesting about these observations, made - no doubt after considerable deliberation - by a theoretical physicist and computer scientist, is that the authority of standard interpretations of classical mathematics is seen, and accepted - perhaps with some element of reluctance, since such acceptance occasionally flies against the grain of observation and experience - not only as absolute, but also as implicitly promising sufficiency, when needed, to help bridge the unbridgeable chasm between a Platonic world of abstract objects, and the real world of sensory perceptions!

## 1.2 Is mathematics the lingua franca of science

As we remark in Anand [An03d], at the heart of this issue is the widespread notion, arguably bordering on misconception, that mathematics is a dispensable tool of science[6], rather than its indispensable mother tongue. If we accept the dictum that what we see reflects more what we are than what there is, then this attitude seems to be reflected, albeit faintly and obliquely, in Shelah's following assertion as well:

> ... some of the best minds in the field of set theory, feel apologetic about their subject. Many are apologetic toward mathematicians, (implying somehow that there are mathematicians and there are logicians, as if they are disjoint species) working in fields which are surely deeper, harder, more profound and meaningful, etc., and so we have to justify our existence by finding applications of "logic" to "mathematics". ... Now, I love to prove theorems in as many areas of mathematics as I can, but I do not like this ... attitude ...

---

[6] An extreme statement of this is seen in Wegner and Goldin's remark that "... computer science is a fundamentally non-mathematical discipline" [WG03].

Now, an explicit thesis of this paper, which is part of a set of related papers where we attempt to view the consequences of Anand [An02a] from a broader perspective, is that such a perception of the nature of mathematics may be less the reflection of any mathematical dogma[7], and more a reflection, as noted in Anand [An03d], of mathematical:

> ... ambiguities, in the classical definitions of foundational elements, that allow the introduction of non-constructive - hence non-verifiable, non-computational, ambiguous, and essentially Platonic - elements into the standard interpretations of classical mathematics[8].

> Consequently, standard interpretations of classical theory may, inadvertently, be weakening a desirable perception - that of mathematics as the lingua franca of scientific expression - by ignoring the possibility that, since mathematics is, indeed, indisputably accepted as the language that most effectively expresses and communicates intuitive truth, the chasm between intuitive truth and formal truth on the one hand, and between formal truth and provability on the other, must, of necessity, be bridgeable.[9]

---

[7] In other words, although standard interpretations of classical mathematical theory may be responsible for tolerating, and perhaps even seeding and sympathetically nurturing, such a perception, they cannot, strictly speaking, be held accountable for its promotion.

[8] See, for instance, the arguments in Anand [An03b].

## 1.3 The roots of Platonism lie in the replacement axiom

That the construction of such a bridge may be a necessity, rather than a matter of philosophical choice, is indicated by Meta-lemma 1 of Anand [An02a], and its consequences, where we essentially prove the following:

> If we give a precise definition of a mathematical object as any individual constant, function letter, or predicate letter that can be introduced into a formal theory without inviting inconsistency, then we can constructively prove that there is a primitive recursive function, $F(x)$, such that the function letter $F$ is not a mathematical object.

An immediate consequence of the above is that the range of such an $F$ cannot be introduced as a set in any axiomatic set theory, which models standard PA, without inviting inconsistency ([An02a], Corollary 1.1). It follows that the replacement axiom of ZFC must, then, be inconsistent with the axioms of standard PA when these are interpreted in ZFC.

Prima facie, accepting the consequences of such a constructively derived result should, reasonably, lead to constructive re-interpretations not only of the obviously non-constructive interpretations of classical set theories, but also to a constructive re-interpretation of the more troublesome classical concepts of Theoretical Computer

---

[9] Of interest, in this context, is Davis's argument [Da95] that an unprovable truth may, indeed, be arrived at algorithmically.

Science that are not obviously non-constructive; for instance, the concept of a constructively defined, but Turing uncomputable, Halting function[10].

More significantly, if Meta-lemma 1 (op cit) implies that we cannot assume that the range of every recursive function defines a (recursively enumerable) set in ZFC ([An02a], Corollary 1.1), then the essential soundness of a set theory based on the standard interpretations of foundational concepts of classical mathematics, where it is implicitly assumed that the range of every recursive function is necessarily a recursively enumerable set in ZFC, may need to be re-assessed. Meta-lemma 1 highlights the fact - which is not obvious, and which does not seem to be explicitly addressed by standard interpretations of classical theory - that such an assumption may be axiomatic; it would, thus, be incapable of proof, although, like the Church-Turing Thesis, it could be vulnerable to disproof in some set theories, as we argue in Anand ([An02a], Corollary 1.1).

Now we note that the replacement axiom was originally intended to prevent the logical, and mathematical, paradoxes from being introduced into axiomatic set theories. However, if the mathematical concepts defined by the antinomies are not mathematical objects (which, in fact, is what the paradoxes establish) then the axiom can possibly be discarded, using Occam's dictum, without any loss of generality.

We also note that, although the replacement axiom may ensure that all the known, non-constructively defined, paradoxical concepts - such as Russell's impredicative set - cannot be introduced as sets of ZFC, it does not, by itself, guarantee that some unknown, constructively defined set, which it may admit as a ZFC set, will not invite inconsistency.

---

[10] Cf. Anand [An03d].

Such an assurance is only available by an appeal to the abstractions of a set-theoretic intuition that, prima facie, appears qualitatively different[11] from the arithmetical intuition that is called upon to accept recursive functions and relations, or the standard interpretations of the axioms of a basic formal system of Arithmetic such as standard PA. However, the arguments of Anand [An02a] imply that our set-theoretic intuition may not be as reliable.

We note, further, that, appearances to the contrary, it is, arguably, the replacement axiom that seems to unrestrictedly admit mathematical concepts that are not mathematical objects into ZFC, rather than the power set axiom, the axiom of choice, or any axiom admitting the existence of a specifically defined infinite set. Thus, it can be reasonably argued that the roots of Platonism may essentially lie in the acceptance of a replacement axiom; if so, then, over time, such of these roots as are critically dependent on the replacement axiom for their sustenance should, therefore, wither away if the axiom is excised.

### 1.4  Standard PA and individual independence / consistency results

We note that accepting, firstly, the introduction of constructive definitions of classical foundational concepts as below, and, secondly, the arguments of Anand [An02a] - for eliminating the replacement axiom - may have some intriguing consequences.

Firstly, there is the possibility that all formal mathematical systems could be constructively built upon standard PA without any loss of generality.

Secondly, assertions in any set-theoretic model of standard PA would necessarily be conditional on the provision of constructive[12] independence, or, at the least, consistency, results for ensuring that the individual constants, function letters, or predicate letters, which occur in such assertions, are, indeed, consistent with the interpretation of the axioms of standard PA in the model (and can thus be introduced into standard PA through appropriate defining axioms[13] without inviting inconsistency).

---

[11] Shelah's remark, "I also have a keen interest in the natural numbers, (though too platonic) but not as a set theorist", also seems to suggest such a qualitative difference.

[12] In contrast to the essentially non-constructive nature of classical forcing arguments.

## 1.5 Shelah's shift to formal proofs of independence

Thus, Shelah's shift from the study of set-theoretic proofs towards the development of formal independence results using forcing methods may be more significant than is apparent from some of his remarks. As we argue in Anand [An02a], not every well-defined mathematical concept corresponds to a mathematical object. Hence, a proof that a concept is consistent with a theory, and can thus be introduced into the theory by suitable defining axioms (that need not be independent) without inviting inconsistency, or a proof that the concept is independent of the theory, and can, therefore, be added to the theory in the form of suitable defining axioms (or their negations) to enlarge the theory without inviting inconsistency, may be crucial for legitimising theorems that are conditional on the concept being a mathematical object.

Shelah's following remarks, regarding his explicit perception of the significance of independence results, are also revealing:

L. Harrington asked me a few years ago: what good does it do you to know all those independences? My answer was: to sort out possible theorems - after throwing away all relations which do not hold you no longer have a heap of questions which clearly are all independent, the trash is thrown away and in what remains you find some grains of gold. This is in general a good justification for independence results; a good place where this had worked is cardinal arithmetic - before Cohen and Easton, who would have looked at $aleph_{omega\_1}^{aleph\_1}$?

---

[13] Cf. Mendelson ([Me64], p82-84).

That this shift goes against the grain of Shelah's instinctive inclination towards set-theoretic methods of proof is reflected in his self-dialogue:

"Does this mean you are a formalist in spite of earlier indications that you are Platonist?" I am in my heart a card-carrying Platonist seeing before my eyes the universe of sets, but I cannot discard the independence phenomena.

## 2. Non-constructivity in classical theory

We thus conjecture that there may have been, and perhaps yet may be, an implicit vision, beyond the horizons of instinctive intuition, that compelled Shaleh to treat independence phenomena as significant; a vision, moreover, which may need to be recognised as a possible beacon. Despite his personal views on the relevance of foundational issues to the study of set theory, Shelah cannot have been completely unaware that there are, indeed, disquieting aspects to the definition of foundational concepts that also cannot be lightly discarded, or ignored.

For instance, in his 1990 article "Second Thoughts About Church's Thesis and Mathematical Proofs", Mendelson's[14] remarks (*italicised parenthetical qualifications added*) suggest that classical definitions of various foundational elements can be argued as being either ambiguous, or non-constructive, or, possibly, both:

Here is the main conclusion I wish to draw: it is completely unwarranted to say that CT (*Church's Thesis*) is unprovable just because it states an equivalence between a vague, imprecise notion (effectively computable function) and a precise mathematical notion (partial-recursive function). ... The concepts and assumptions that support the notion of partial-recursive function are, in an essential way, no less vague and imprecise (*non-constructive, and intuitionistically objectionable*) than the notion of effectively computable function; the former are just more familiar and are part of a respectable theory with connections to other parts of logic and mathematics. (The notion of effectively computable function could have been incorporated into an axiomatic presentation of classical mathematics, but the acceptance of CT made this unnecessary.) ... Functions are defined in terms of sets, but the concept of set is no clearer (*not more non-constructive, and intuitionistically objectionable*) than that of function and a foundation of mathematics can be based on a theory using function as primitive notion instead of set. Tarski's definition of truth is formulated in set-theoretic terms, but the notion of set is no clearer (*not more non-constructive, and intuitionistically objectionable*), than that of truth. The model-theoretic definition of logical validity is based ultimately on set theory, the foundations of which are no clearer (*not more non-constructive, and intuitionistically objectionable*) than our intuitive (*non-constructive, and intuitionistically objectionable*) understanding of logical validity. ... The notion of Turing-computable function is no clearer (*not more non-constructive, and intuitionistically objectionable*) than, nor more mathematically useful (foundationally speaking) than, the notion of an effectively computable function.

---

[14] Mendelson [Me90].

## 2.1 Constructive definitions

In Anand [An02a], we argue that accepting such leeway in the standard interpretations of classical concepts leads to various ambiguities and anomalies that may, however, be avoidable if we define foundational concepts - such as the following - constructively in terms of a smaller number of primitive, formally undefined but intuitively unobjectionable, mathematical concepts:

(*i*) **Primitive mathematical object**: A primitive mathematical object is any symbol for an individual constant, predicate letter, or a function letter, which is defined as a primitive symbol of a formal mathematical language.[15]

(*ii*) **Formal mathematical object**: A formal mathematical object is any symbol for an individual constant, predicate letter, or a function letter that is either a primitive mathematical object, or that can be introduced through definition into a formal mathematical language without inviting inconsistency.[16]

(*iii*) **Mathematical object**: A mathematical object is any symbol that is either a primitive mathematical object, or a formal mathematical object.

(*iv*) **Set**: A set is the range of any function whose function letter is a mathematical object.

(*v*) **Individual computability**: A number-theoretic function $F(x)$ is individually computable if, and only if, given any natural number $k$, there is an individually effective method (which may depend on the value $k$) to compute $F(k)$.

(*vi*) **Uniform computability**: A number-theoretic function $F(x)$ is uniformly computable if, and only if, there is a uniformly effective method (necessarily independent of $x$) such that, given any natural number $k$, it can compute $F(k)$.

(*vii*) **Effective computability**: A number-theoretic function is effectively computable if, and only if, it is either individually computable, or it is uniformly computable.[17]

---

[15] We note that, as remarked by Mendelson [Me90], the terms "function" and "function letter" - and, presumably, "individual constant", "predicate", and "predicate letter" - can be taken as undefined, primitive foundational concepts.

[16] We highlight the significance of this definition in Meta-lemma 1 in Anand [An02a].

[17] We note that classical definitions of the effective computability of a function (cf. [Me64], p207) do not distinguish between the two cases. The standard interpretation of effective computability is to implicitly

(*viii*) **Individual truth**: A string $[F(x)]^{18}$ of a formal system P is individually true under an interpretation M of P if, and only if, given any value $k$ in M, there is an individually effective method (which may depend on the value $k$) to determine that the interpreted proposition $F(k)$ is satisfied in M.[19]

(*ix*) **Uniform truth**: A string $[F(x)]$ of a formal system P is uniformly true under an interpretation M of P if, and only if, there is a uniformly effective method (necessarily independent of $x$) such that, given any value $k$ in M, it can determine that the interpreted proposition $F(k)$ is satisfied in M.

(*x*) **Effective truth**: A string $[F(x)]$ of a formal system P is effectively true under an interpretation M of P if, and only if, it is either individually true in M, or it is uniformly true in M.[20]

(*xi*) **Individual Church Thesis**: If, for a given relation $R(x)$, and any element $k$ in some interpretation M of a formal system P, there is an individually effective method such that it will determine whether $R(k)$ holds in M or not, then every element of the domain D of M is the interpretation of some term of P, and there is some P-formula $[R'(x)]$ such that:

---

treat it as equivalent to the assertion: A number-theoretic function $F(x)]$ of a formal system P is effectively computable if, and only if, it is both individually computable, and uniformly computable.

[18] We use square brackets to distinguish between the uninterpreted string $[F]$ of a formal system, and the symbolic expression "$F$" that corresponds to it under a given interpretation that unambiguously assigns formal, or intuitive, meanings to each individual symbol of the expression "$F$".

[19] In Anand [An02a], we argue that, under a constructive interpretation of formal Peano Arithmetic, Gödel's undecidable proposition, is individually, but not uniformly, true under the standard interpretation. See also Anand [An03c].

[20] We note that, classically, Tarski's definition of the truth of a formal proposition under an interpretation (cf. [Me64], p49-52) does not distinguish between the two cases. The implicitly accepted (standard) interpretation of the definition appears, prima facie, to be the non-constructive assertion: A string $[F(x)]$ of a formal system P is true under an interpretation M of P if, and only if, it is both uniformly true in M, and individually true in M.

$R(k)$ holds in M if, and only if, $[R'(k)]$ is P-provable.

(In other words, the Individual Church Thesis postulates that, if a relation $R$ is effectively decidable individually (possibly non-algorithmically) in an interpretation M of some formal system P, then $R$ is expressible in P, and its domain necessarily consists of only mathematical objects, even if the predicate letter $R$ is not, itself, a mathematical object.)

(*xii*) **Uniform Church Thesis**: If, in some interpretation M of a formal system P, there is a uniformly effective method such that, for a given relation $R(x)$, and any element $k$ in M, it will determine whether $R(k)$ holds in M or not, then $R(x)$ is the interpretation in M of a P-formula $[R(x)]$, and:

$R(k)$ holds in M if, and only if, $[R(k)]$ is P-provable.

(Thus, the Uniform Church Thesis postulates that, if a relation $R$ is effectively decidable uniformly (necessarily algorithmically) in an interpretation M of a formal system P, then, firstly, $R$ is expressible in P, and, secondly, the predicate letter $R$, and all the elements in the domain of the relation $R$, are necessarily mathematical objects.)

## 2.2 Some consequences

The significance of constructively interpreting foundational concepts and assertions of classical mathematics is that:

(*i*) The Uniform Church Thesis implies that a formula $[R]$ is P-provable if, and only if, $[R]$ is uniformly true in some interpretation M of P.

(*ii*) The Uniform Church Thesis implies that if a number-theoretic relation $R(x)$ is uniformly satisfied in some interpretation M of P, then the predicate letter "$R$" is a formal mathematical object in P (i.e. it can be introduced through definition into P without inviting inconsistency).

(*iii*) The Uniform Church Thesis implies that, if a P-formula [*R*] is uniformly true in some interpretation M of P, then [*R*] is uniformly true in every model of P.

(*iv*) The Uniform Church Thesis implies that if a formula [*R*] is not P-provable, but [*R*] is classically true under the standard interpretation, then [*R*] is individually true, but not uniformly true, in the standard model of P.

(*v*) The Uniform Church Thesis implies that Gödel's undecidable sentence GUS is individually true, but not uniformly true, in the standard model of P.[21]

By defining effective computability, both individually and uniformly, along similar lines, we can give a constructive definition of uncomputable number-theoretic functions:

(*vi*) A number-theoretic function $F(x_1, ..., x_n)$ in the standard interpretation M of P is uncomputable if, and only if, it is effectively computable individually, but not effectively computable uniformly.

This, last, removes the mysticism behind the fact that we can define a number-theoretic Halting function that is, paradoxically, Turing-uncomputable.

(*vii*) If we assume a Uniform Church Thesis, then every partial recursive number-theoretic function $F(x_1, ..., x_n)$ has a unique constructive extension as a total function.

(*viii*) If we assume a Uniform Church Thesis, then not every effectively computable function is classically Turing computable (so Turing's Thesis does not, then, hold).

---

[21] An intriguing consequence of this argument is considered in Appendix 1 of Anand [An03c].

(*ix*) If we assume a Uniform Church Thesis, then not every (partially) recursive function is classically Turing-computable.[22]

(*x*) If we assume a Uniform Church Thesis, then the class P of polynomial-time languages in the P versus NP problem may not define a formal mathematical object.

## 3. Generality and the Entscheidungsproblem

An intriguing aspect of Shelah's lecture is his revealing remark, reproduced below, that, despite his shift towards the development of individual independence results, his ideological preference is yet for general methods of proof.

Now, Hibert's Entscheidungsproblem can, arguably, be seen as heralding a paradigm shift towards an accelerated focusing on the study of general methods of proof - the traditional prerogative of pure mathematicians - at the possible expense of a reduced attention to the development of specific solutions for individual cases. The latter are, prima facie, the primary concern of the applied sciences, which rely on mathematics to supply a suitably verifiable language for expressing the results of their specific observations, albeit in a coherent and broadly predictable form

Ironically, the negative answer to the Entscheidungsproblem, which emerged out of the work of Church, Turing and Goedel, can, again arguably, be seen as having been perceived as a challenge by pure mathematics to concentrate even more on staking out, and fortressing, those areas where general expressions can be meaningfully asserted. An arguable consequence of such a siege mentality: the interpretation of concepts in terms that appear, momentarily, to have an essentially non-existent relatability, in some cases[23],
to the original concerns that prompted the consideration of such ideas in the first place.

---

[22] The classical proof that every (partially) recursive function is classically Turing-computable uses induction over (partial) recursive functions, thus assuming that every such function is a mathematical object; by Meta-lemma 1, such an assumption is invalid.

That a bias towards the consideration of general methods may be instinctive, in some
cases, is suggested by Shelah's remarks:

> I was attracted first to mathematics and subsequently to mathematical logic by their
> generality, anticipating that this is the normal attitude; it seems I was mistaken. I have
> always felt that examples usually just confuse you (though not always), having
> always specific properties that are traps as they do not hold in general. Note that by
> "generality" I mean I prefer, e.g., to look more at general complete first order theories
> (possibly uncountable) rather than at simple groups of finite Morley rank.

> However I do not believe in "never look at the points, always look at the arrows";
> each problem has to be dealt with according to its peculiarities, and finding
> applications of your own field in another means showing something that interests the
> others; but given a problem, why not try for the best, most general statement available
> (of course, if the theorem exists, and the additional generality requires no substance, it
> is not exciting).

Awareness of an instinctive preference for general methods is also seen in Shelah's
response to being "accused" of:

> ... (*explaining*) in detail why proofs in ZFC are best, and why I prefer them to
> independence results, just two years before launching full scale into forcing.

where he notes that:

---

[23] E.g. the admittedly non-intuitive notion of "inaccessible" cardinals; at least until such time that they can
be made "accessibl" in some sense by a constructive definition of the notion.

I still feel an outright answer in ZFC is best, even though a new technique for proving independence may be more interesting. Cohen's theorem seems to me more interesting than a proof of CH, as it supplies us with a general method.

## 3.1 Has this bias adversely affected the development of mathematics?

Now, it can be argued that such a bias, towards an overwhelming development of general methods of proof, may be seriously inhibiting the development of mathematics as a universal language, which can effectively express all forms of human cognition.

Increasingly, the major challenges of Theoretical Computer Science, Quantum Physics, and other disciplines are, prima facie, to express what appear to be non-algorithmic, but determinate, processes; such processes seem to characterise natural laws more than classical algorithmic processes. However, the implicit thesis of standard interpretations of classical theory - reflected in the broad acceptance of CT - that such phenomena are essentially inaccessible to effective methods of expression seems, prima facie, to go against the growing body of experimental evidence to the contrary[24].

As we remark in Anand [An03c]:

... the central issue in the development of AI is that of finding effective methods of duplicating the cognitive and expressive processes of the human mind. This issue is being increasingly brought into sharper focus by the rapid advances in the experimental, behavioral, and computer sciences[25]. Penrose's "The Emperor's New Mind", and "Shadows of the Mind", highlight what is striking about the attempts, and struggles, of current work in these areas to express their observations adequately - necessarily in a predictable way - within the standard interpretations of formal propositions as offered by classical theory.

---

[24] In his 2003 BBC Reith lectures, Ramachandran speculates that it may be possible, at some future date, to map, and physically link, the cognitive parts of one brain into another, so the latter can mirror the former's sensory perceptions as identically sensed experiences.

[25] See, for instance, footnote 18 in Ramachandran.[RH01].

So, the question arises: Are formal classical theories essentially unable to adequately express the extent and range of human cognition, or does the problem lie in the way formal theories are classically interpreted at the moment? The former addresses the question of whether there are absolute limits on our capacity to express human cognition unambiguously; the latter, whether there are only temporal limits - not necessarily absolute - to the capacity of classical interpretations to communicate unambiguously that which we intended to capture within our formal expression.

The thesis of this, and related, papers[26] is that we may comfortably reject the former by recognising, firstly, that we can, indeed, constructively define foundational concepts unambiguously as indicated above, and, secondly, that, appearances to the contrary, all set-theoretic concepts should be capable of constructive interpretations without any loss of generality.

---

[26] See Anand, in particular [An02a], [An02b], [An03a], [An03b], [An03c], and [An03d].

(*Updated: Saturday 10[th] May 2003 8:30:28 AM IST by re@alixcomsi.com*)